\documentclass{article}
\usepackage[utf8]{inputenc}
\usepackage{amssymb}
\usepackage{amsmath}
\usepackage{amsthm}
\usepackage{color}
\usepackage{enumitem}
\usepackage{hyperref}
\usepackage[margin=1.3in]{geometry}
\usepackage{tikz}
\usepackage{bbm}

\theoremstyle{definition}
\newtheorem{definition}{Definition}[section]

\newtheorem{remark}[definition]{Remark}
\theoremstyle{plain}
\newtheorem{thm}[definition]{Theorem}

\newtheorem{prop}[definition]{Proposition}
\newtheorem{lemma}[definition]{Lemma}

\newtheorem{coro}[definition]{Corollary}

\newcommand{\vol}{\text{Vol}}
\newcommand{\sys}{\text{sys}}

\title{Systolic inequalities for $\mathbb{S}^1$-invariant contact forms in dimension three}
\author{Simon Vialaret\thanks{Université Paris-Saclay, CNRS, Laboratoire de mathématiques d'Orsay, 91405, Orsay, France \\simon.vialaret@universite-paris-saclay.fr}}
\date{}

\begin{document}

\maketitle

\begin{abstract}
	In contact geometry, a systolic inequality is a uniform upper bound on the shortest period of a closed Reeb orbit, in terms of the contact volume. We prove a general systolic inequality valid on Seifert bundles with non-zero Euler number for all contact forms that are invariant under the underlying circle action.
\end{abstract}

\section{Introduction}

\subsection{Systolic inequalities in contact geometry} % (fold)
\label{sub:systolic_inequalities_in_contact_geometry}

A \textbf{contact form} on a three-manifold $M$ is a one-form $\alpha$ such that $\alpha \wedge \mathrm{d}\alpha$ is a volume form. A contact form $\alpha$ induces an orientation on $M$ given by $\alpha \wedge \mathrm{d}\alpha$. It also induces a \textbf{contact structure}, meaning a nowhere integrable plane distribution, given by $\xi = \ker \alpha$. Finally, the \textbf{Reeb vector field} of $\alpha$, written $R$, is uniquely defined by the equations $\iota_R \alpha =1$ and $\iota_R \mathrm{d}\alpha=0$. The flow of $R$ is the \textbf{Reeb flow} of $\alpha$. The Weinstein conjecture states that all Reeb flows on a closed manifold have a periodic orbit. In dimension three, it was proved by Taubes in \cite{10.2140/gt.2007.11.2117}. This fundamental fact allows to define the \textbf{systole} of a contact form $\alpha$ on a three-manifold $M$, that we denote $\sys(\alpha)$, as the infimum of the set of periods of periodic orbits of the Reeb flow of $\alpha$. The \textbf{systolic ratio} of $\alpha$, written $\rho(\alpha)$, is the scaling-invariant quantity $\rho(\alpha)=\frac{\sys(\alpha)^2}{\vol(\alpha)}$, where $\vol(\alpha) = \int_M \alpha \wedge \mathrm{d}\alpha$ is the contact volume of $M$. Finally, a subset $\Omega$ of the set of contact forms on $M$ is said to satisfy a \textbf{systolic inequality} whenever the systolic ratio is bounded on $\Omega$, in other words if there is a positive constant $C_\Omega$ such that for all contact form $\alpha \in \Omega$, the following inequality holds:
\[
	\sys(\alpha)^2 \leq C_\Omega \vol(\alpha).
\]
	
Systolic inequalities for contact forms have interesting applications. First, a rich class of examples of contact forms is provided by Riemannian and Finsler metrics, since a metric on a manifold $B$ induces a canonical contact form on the unit tangent bundle of $B$, whose Reeb flow coincides with the geodesic flow. Systolic inequalities have been widely studied in Riemannian and Finsler geometry, see for example \cite{Katz_2007} and the survey \cite{benedetti2021steps}. Since \cite{paiva_2014_contact}, techniques from contact geometry have been introduced in systolic geometry, and have shed some new light on some problems, see \cite{paiva_2014_contact}, \cite{Abbondandolo2017}, \cite{Abbondandolo2018}, \cite{ABHS_revol}, \cite{benedetti_kang_2020}, \cite{abbondandolo:ensl-03357629}, \cite{MR4529024}, \cite{vialaret2024sharpsystolicinequalitiesinvariant} for some recent examples.

A second motivation for the study of systolic inequalities in contact geometry is their relation with some symplectic capacities. Symplectic capacities are maps sending open subsets of $\mathbb{R}^{2n}$ to $[0, + \infty]$ which are invariant under symplectomorphisms, monotone under inclusion, homogeneous of degree 2 with respect to dilations, and satisfying a normalization axiom. In \cite{Viterbo_2000}, Viterbo states the following conjecture: for any symplectic capacity $c$, and for any bounded convex domain $U \subset \mathbb{R}^{2n}$, the inequality
\begin{equation}
	\label{eq:viterbo_conjecture}
	c^n(U) \leq n! \vol(U)
\end{equation}
holds. Moreover, equality holds if and only if the interior of $K$ is symplectomorphic to a ball. This conjecture has been proved to be wrong at this level of generality by Haim-Kislev and Ostrover in \cite{haimkislev2024counterexampleviterbosconjecture}, in which they construct examples of convex domains in any dimension larger than or equal to 4 whose Hofer--Zehnder capacity fails to satisfy Viterbo's conjecture. However, some known consequences of this conjecture, such as a conjecture of Mahler in convex geometry (see \cite{10.1215/00127094-2794999}), are still open and are strong motivations to study its domain of validity.

If a convex domain $U \subset \mathbb{R}^{2n}$ has smooth boundary, then $\partial U$ has contact type. As such, it carries a natural contact form $\alpha_U$. It is a striking fact that many known capacities coincide with the systole of $\alpha_U$ when evaluated on a convex domain $U$ (see \cite{hryniewicz2023hopf} and \cite{abbondandolo2024closedcharacteristicsminimalaction} for some precise statements). For any such capacity, inequality (\ref{eq:viterbo_conjecture}) is equivalent to the systolic inequality $\sys(\alpha_U)^{n} \leq \vol(\alpha_U)$ for any convex domain $U$.

It is a striking fact that the class of all contact forms on any given compact manifold $M$, adapted to any given contact structure $\xi$ does not satisfy a systolic inequality. More precisely, for any $\epsilon > 0$, there exists a contact form $\alpha$ on $M$ such that
\begin{itemize}
	\item $\alpha$ induces the contact structure $\xi$;
	\item $\alpha$ has contact volume less than $\epsilon$;
	\item $\alpha$ has systole larger than one.
\end{itemize}
This was proved in \cite{Abbondandolo2018} for the standard contact structure on the three-sphere, in \cite{Abbondandolo2019} in dimension three, and in \cite{Saglam_2021} in higher dimensions. The goal of this work is to show that the presence of sufficient symmetries of the Reeb flows forces the existence of systolic inequalities.

\subsection{Seifert bundles, invariant contact forms and the main result} % (fold)
\label{sub:invariant_contact_forms}

In this work, we consider the class of contact forms which are invariant under an almost-free circle action on a three-manifold $M$, namely a circle action on $M$ such that all orbits have trivial stabilizers, except for finitely many with finite stabilizer. Such an action gives $M$ the structure of a Seifert bundle. The classification achieved in \cite{10.1007/BF02398271} states that Seifert bundles can be obtained by performing a finite number of Dehn surgeries on unknots in $\mathbb{S}^1 \times S_g$, where $S_g$ is a closed oriented surface with genus $g$. Each surgery is described by coprime coefficients $(p_i,q_i) \in \mathbb{N} \times \mathbb{Z}$. Furthermore, whether two data $(g, (p_i, q_i)_{1 \leq i \leq L})$ and $(g', (p'_i, q'_i)_{1 \leq i \leq L'})$ describe the same Seifert bundle is fully understood and is recalled in Theorem \ref{thm:classification_seifert} below. In particular, the \textbf{Euler number} $e$ of a Seifert bundle, defined by $e = -\sum_{l=1}^L \frac{q_l}{p_l}$, is a well-defined invariant, which vanishes if and only if $M$ is finitely covered by a trivial circle bundle. In what follows, $M(g, (p_i, q_i)_{1 \leq i \leq L})$ will stand for the Seifert bundle obtained after performing the surgeries with coefficients $(p_i, q_i)_{1 \leq i \leq L}$ on $\mathbb{S}^1 \times S_g$. We will recall some basic facts about Seifert bundles in Section \ref{sub:seifert_bundles}.

Writing $X$ for the vector field generating the circle action on $M$, we say that a contact form $\alpha$ is \textbf{invariant} whenever $\mathcal L_X \alpha = 0$. We write $\Omega(g, (p_i, q_i)_{1 \leq i \leq L})$ for the class of invariant contact forms on $M(g, (p_i, q_i)_{1 \leq i \leq L})$.

A source of examples of invariant contact forms is given by the class of \textbf{Besse} contact forms, which are the contact forms all of whose Reeb orbits are closed. The Reeb flow of a Besse form $\alpha_0$ induces an almost-free circle action on the ambient manifold, which leaves $\alpha_0$ invariant. When this action is free, namely when all Reeb orbits have the same minimal period, the contact form is said to be \textbf{Zoll}. Zoll contact forms have a special role in the study of the systolic ratio. Indeed, a result of Àlvarez-Paiva and Balacheff \cite{paiva_2014_contact} states that a contact form that is a local maximizer of the systolic ratio is necessarily Zoll. Conversely, it was proved in \cite{Abbondandolo_Benedetti_2023} that all Zoll contact forms are local maximizers of the systolic ratio. A similar statement for Besse contact forms is also true, where the systolic ratio is replaced by higher systolic ratios, see Theorem A in \cite{abbondandolo:ensl-03357629} for a precise statement, and \cite{baracco2023local} for a similar result for the higher Ekeland--Hofer capacities.

For Besse contact forms, the circle action is everywhere transverse to the contact structure. This does not need to be the case for a general contact structure invariant under a $\mathbb{S}^1$-action. In fact, $\Omega(g, (p_i, q_i)_{1 \leq i \leq L})$ has infinitely many connected components, and for most of them, the circle action has a non-empty collection of orbits that are Legendrian. More details are given in Section \ref{sub:classification_of_invariant_contact_structures}, where the classification of $\mathbb{S}^1$-invariant contact structures from \cite{niederkruger2005} in terms of the collection of Legendrian orbits of the action is recalled without a proof.

In \cite{vialaret2024sharpsystolicinequalitiesinvariant}, the author proved a sharp systolic inequality for invariant tight contact forms on principal circle bundles over $\mathbb{S}^2$, generalizing the main result of \cite{ABHS_revol} about Riemannian metrics of revolution on $\mathbb{S}^2$. The main result of this work is a general systolic inequality holding for invariant contact forms on Seifert bundles with non-vanishing Euler number.

\begin{thm}
	\label{thm:ineg_sys}
	There exists a positive constant $C > 0$ such that the following holds:

	For all $g \in \mathbb{N}$ and all finite families of coprime integers $(p_i, q_i)_{1 \leq i \leq L}$ satisfying $e = -\sum \limits_{i=1}^L \frac{q_i}{p_i} \neq 0$, we have
	\begin{equation}
		\label{eq:ineg_sys}
		\forall \alpha \in \Omega(g,(p_i, q_i)_{1 \leq i \leq L}), \hspace{0.5cm} \sys(\alpha)^2 \leq C\max\left(1, \frac{1}{|e|}\right) \vol(\alpha).
	\end{equation}
\end{thm}

\begin{remark}

	\hspace{0.01cm}
	\begin{itemize}
		\item Let us comment on the aymptotic behavior of the systolic inequality (\ref{eq:ineg_sys}) when $e$ goes to $+\infty$ and to zero. First, it was proved in \cite[Proposition 4.1]{vialaret2024sharpsystolicinequalitiesinvariant} that for all integers $e > 2$, there exists an invariant contact form on the total space of the principal bundle over $\mathbb{S}^2$ of Euler number $e$, with systolic ratio larger than $\frac{1}{4}$. It follows that the right-hand side in (\ref{eq:ineg_sys}) has the right asymptotic behavior when $e$ goes to $+ \infty$, since it is constant. However, it is not clear whether there is a sequence of Seifert bundles $(M_n)_{n \in \mathbb{N}}$ with Euler numbers $(e_n)_{n \in \mathbb{N}}$ going to zero, each of them carrying an invariant contact form $\alpha_n$, such that $(\rho(\alpha_n))_{n \in \mathbb{N}}$ goes to $+\infty$, as suggested by the right-hand side of \ref{eq:ineg_sys}.
		\item Although the fact that the Euler number is non-zero is a key argument in the proof of Theorem \ref{thm:ineg_sys}, it is still unclear whether or not a systolic inequality holds for invariant contact forms when $e = 0$. For example, it is not clear how to construct invariant contact forms with arbitrarily large systolic ratio on $\mathbb{S}^1 \times \mathbb{S}^2$.
		\item The assumption of being invariant under an almost-free circle action is essential for the statement to be true. Indeed, given any integer $k$, one can use the methods of \cite{Abbondandolo2018} to produce contact forms on $\mathbb{S}^3$ with arbitrarily large systolic ratio that are invariant under a free $\mathbb{Z}/k \mathbb{Z}$-action. Moreover, those contact forms can be chosen so that their Reeb flows are Liouville integrable.
		\item The constant $C$ can be chosen to be $225$. This is certainly not the optimal constant. A sharp systolic inequality for $\mathbb{S}^1$-invariant contact forms was proved under more restrictive conditions by the author in \cite{vialaret2024sharpsystolicinequalitiesinvariant}.
	\end{itemize}
\end{remark}

This theorem has as an immediate corollary a systolic inequality for invariant contact forms on non-trivial principal bundles, as in this case $e$ is an integer.

\begin{coro}
	There exists a positive constant $C > 0$ such that for all invariant contact form $\alpha$ on the total space of a non-trivial $\mathbb{S}^1$-principal bundle over a closed orientable surface,
	\[
		\sys(\alpha)^2 \leq C \vol(\alpha).
	\]
\end{coro}

% subsection context_and_main_result (end)

\subsection{A sketch of proof} % (fold)
\label{sub:a_sketch_of_proof}

We now present a sketch of proof of Theorem \ref{thm:ineg_sys}. First of all, the invariance of the Reeb flow by a group action allows to produce a non-trivial first integral of the dynamic, namely an equivariant function $K \colon M \rightarrow \mathbb{R}$ that is preserved by the Reeb flow. Whenever this first integral is nowhere vanishing, an elementary argument yields a sharp systolic inequality. This is the content of Proposition \ref{prop:ineg_sys_facile}. When $K$ is not nowhere-vanishing, more ingredients are needed.

First, we construct a family of equivariant surfaces of section for Reeb flows of invariant contact forms. A surface of section for a flow $(\phi_t)_{t \in \mathbb{R}}$ in a three-manifold is an embedded surface $\Sigma$ transverse to the flow and which has the following property: for all $x \in \Sigma$, there exists a positive time $\tau(x)$ and a negative time $\bar \tau (x)$ such that $\phi_{\tau(x)}(x)$ and $\phi_{\bar\tau(x)}(x)$ belong to $\Sigma$. While the existence of a surface of section for a general flow is a difficult problem, they are a powerful tool for the study of three-dimensional dynamics which allows to reduce the problem to the well-studied setting of surface diffeomorphisms. The search for surfaces of section for Reeb flows has been a driving force in contact dynamics, and has had many fruitful applications, in particular to the study of the systolic ratio, see the survey \cite{SALOMAO_HRYNIEWICZ_2019} and \cite{Abbondandolo2018}, \cite{Saglam_2021}, \cite{ABHS_revol}, \cite{benedetti_kang_2020}, \cite{vialaret2024sharpsystolicinequalitiesinvariant} for some recent applications.

We prove in Section \ref{sub:invariant_surface_of_section} that Reeb flows of invariant contact forms admit a surface of section $\Sigma$ which is a finite disjoint union of annuli, and which is furthermore invariant under the circle action. That surface of section has as many components as the union of the Legendrian $\mathbb{S}^1$-orbits of the underlying invariant contact structure. Then, in Section \ref{sub:a_family_of_potentials}, we use this invariant surface of section to reduce the study of most of the Reeb flow to a one-dimensional problem. More precisely, for each component of $\Sigma$, we construct a real-valued function, defined on an interval, which we call a \textbf{potential}, and which encodes the first-return map of that component. The precise dictionary between properties of the first-return map, in particular its periodic points and their action, and properties of the potential, is stated in Corollary \ref{coro:periodic_point}. The construction of the potentials makes use of a specific decomposition of the contact form which we prove in Section \ref{sub:a_decomposition_for_invariant_contact_forms}.

We can now proceed by contradiction, and assume that there is no systolic inequality on $\Omega(g, (p_i, q_i)_{1 \leq i \leq L})$. By definition, this means that there are invariant contact forms whose systolic ratio is arbitrarily large. The previous paragraph tells us that we can associate to any such contact form a family of potentials. The property of having high systolic ratio can be translated to strong constraints on the potentials in terms of the variations of their derivatives, and of their integrals via Corollary \ref{coro:periodic_point} and Lemma \ref{lemma:volume_J}. We prove in Section \ref{sec:necessary_properties_high_systolic_ratio} that those constraints ultimately imply that each of the potentials must be $\mathcal C^1$-close to linear maps, meaning that the first-return maps of each component of $\Sigma$ must be uniformly close to rotations of the annulus. However, we prove in Lemma \ref{lemma:computing_euler} an identity relating the Euler number of the Seifert bundle and the potentials. The contradiction now follows from the fact that the only value for the Euler number that allows the family of potentials to be uniformly close to linear maps is zero. Roughly speaking, for the surfaces of section to glue up to produce a \emph{non-trivial} circle bundle, the first-return maps must be different enough from uniform rotations. Hence, a contradiction, which gives the systolic inequality.

% subsection a_sketch_of_proof (end)
\begin{remark}
	In what follows, any object defined in term of a contact form $\alpha$ will have a superscript $\alpha$. This superscript will be omitted if the context makes clear which contact form is being used. Objects defined on the base of a bundle will carry a tilde (\textasciitilde) symbol, while lifts to the total space will be written without tilde.
\end{remark}

\paragraph{Acknowledgements}
The author is very grateful to Alberto Abbondandolo and Rémi Leclercq for numerous stimulating conversations and advice. The author is partially funded by the Deutsche Forschungsgemeinschaft (DFG, German Research Foundation) – Project-ID 281071066 – TRR 191 and by the ANR grant 21-CE40-0002 (CoSy).

\section{Invariant contact forms on Seifert bundles} % (fold)
\label{sec:invariant_contact_forms_on_seifert_bundles}

In this section, we give more background on Seifert bundles, as well as an overview of the classification of invariant contact structures on such spaces.

\subsection{Seifert bundles, Euler number and normalized invariants} % (fold)
\label{sub:seifert_bundles}

From now on, we will call \textbf{Seifert bundle} a closed oriented manifold of dimension three, endowed with an effective fixed-point free smooth action of the circle. Although this only corresponds to the case of orientable total space and orientable base (``Oo'' in the original notation of H. Seifert) in the classification of Seifert spaces from \cite{10.1007/BF02398271}, we will keep the name Seifert bundle for the sake of brevity. We will denote by $\pi$ the projection from $M$ to the orbit space of the action. Orbits of the action are sorted in two classes, depending on the stabilizer $G_x$ of their base point $x$. If $G_x$ is trivial, the orbit of $x$ under the circle action is called \textbf{regular}. Otherwise, $G_x$ is a finite cyclic group $\mathbb{Z}/p \mathbb{Z}$ and in that case the orbit is called \textbf{singular}. It follows from the effectiveness of the action and from the slice theorem for smooth actions of compact groups that singular orbits are in finite number. Moreover, a singular orbit with stabilizer $\mathbb{Z} / p \mathbb{Z}$ has a neighborhood $V$ which is equivariantly diffeomorphic to the local model $D^2 \times \mathbb{S}^1 \rightarrow D^2$ given by $(rt_1, t_2) \mapsto r t_1^p t_2^q$, where $p$ and $q$ are coprime integers.

An alternative way to describe Seifert bundles is by the successive use of Dehn surgeries on a trivial circle bundle. Let $S_g$ be an oriented closed surface of genus $g$, $\gamma$ a $\mathbb{S}^1$-fiber of $\mathbb{S}^1 \times S_g \overset{\pi}{\rightarrow}S_g$, $U$ a tubular neighborhood of $\gamma$ which is a union of $\mathbb{S}^1$-fibers, and $(p,q) \in \mathbb{Z}^2$ two coprime integers with $p > 0$. The boundary of $M \setminus U$ is a 2-torus whose fundamental group is spanned by $(l,m)$, where $l$ is homotopic to a positively oriented $\mathbb{S}^1$-fiber in $M \setminus \gamma$ and $-\pi_*m$ is a positively oriented boundary component of $\partial (S \setminus \pi(U))$. Similarly let $T = \mathbb{S}^1 \times D^2$ be an oriented solid torus. The fundamental group of $\partial T$ is generated by $(\tilde l,  \tilde m)$, where $\tilde l$ is the positive generator of $\pi_1(T)$ and $\tilde m$ is homotopic to $\{t_0\} \times \partial D^2$ in $\partial T$, with $t_0 \in \mathbb{S}^1$, and oriented as $\partial D^2$. Gluing $T$ to $M \setminus \gamma$ along their boundaries via a diffeomorphism of $\mathbb{T}^2$ identifying $m$ to $p \tilde m + q \tilde l$ produces a Seifert bundle which is a genuine circle bundle when $p=1$, and which has a singular orbit with stabilizer $\mathbb{Z}/ p \mathbb{Z}$ otherwise. Iterating this process allows us to define $M(g,(p_i, q_i)_{1 \leq i \leq L})$ as the Seifert bundle obtained by doing a finite number of Dehn surgeries on regular orbits with gluing coefficients $(p_i,q_i)_{1\leq i \leq L}$. It is not hard to prove that any Seifert bundle can be obtained by doing so. The classification in \cite{10.1007/BF02398271} (see also the lecture notes \cite{jankins1983lectures}) describes when two collections of coefficients $(p_i,q_i)$ produce the same Seifert bundle.

\begin{thm}[Theorem 1.5 in \cite{jankins1983lectures}]
	\label{thm:classification_seifert}
	$M(g,(p_i, q_i)_{1 \leq i \leq L})$ and $M(g',(p'_i,q'_i)_{1 \leq i \leq L'})$ are equivariantly diffeomorphic if and only if the following conditions are satisfied:
	\begin{enumerate}
		\item $g=g'$;
		\item disregarding any $\frac{q_l}{p_l}$ and $\frac{q'_l}{p'_l}$ which are integers, the remaining $\frac{q_l}{p_l} (\text{mod }1)$ are a permutation of the remaining $\frac{q'_l}{p'_l} (\text{mod }1)$;
		\item $\sum_{l=1}^L \frac{q_l}{p_l} = \sum_{l=1}^{L'} \frac{q'_l}{p'_l}$.
	\end{enumerate}
\end{thm}

The \textbf{Euler number} of $M(g,(p_i, q_i)_{1 \leq i \leq L})$ is defined by $e = -\sum_{l=1}^L \frac{q_l}{p_l}$. It follows from Theorem \ref{thm:classification_seifert} that is it an invariant of Seifert bundles, since it does not depend on the choice of a surgery description. When the bundle does not have singular orbits, then $e$ is the usual Euler number of the corresponding $\mathbb{S}^1$-principal bundle.

We will say that a family $(p_i, q_i)_{1 \leq i \leq L}$ is \textbf{normalized} if $L=1$ in the case where there is no singular orbits, and if $L$ equals to the number of singular orbits otherwise. The classification in Theorem \ref{thm:classification_seifert} implies that given $(p_i, q_i)_{1 \leq i \leq L}$ not necessarily normalized, we can always find $(p'_i, q'_i)_{1 \leq i \leq L'}$ normalized such that $M(g,(p_i, q_i)_{1 \leq i \leq L})$ and $M(g,(p'_i, q'_i)_{1 \leq i \leq L'})$ are equivariantly diffeomorphic. From now on we will always assume $(p_i, q_i)_{1 \leq i \leq L}$ to be normalized.

\subsection{Invariant contact forms and first properties} % (fold)
\label{sub:invariant_contact_forms_and_first_properties}

Let $M = M(g,(p_i, q_i)_{1 \leq i \leq L})$ be a Seifert manifold, and $X$ the vector field generating the circle action on $M$. We recall that a contact form $\alpha$ on $M$ is \textbf{invariant} if it coincides with its pull-back by the circle action, which is equivalent to $\mathcal L_X \alpha=0$. We denote by $\Omega(g,(p_i, q_i)_{1 \leq i \leq L})$ the set of invariant contact forms on $M(g,(p_i, q_i)_{1 \leq i \leq L})$. It follows from the invariance that the Reeb flow of $\alpha$ preserves the moment map $K^\alpha = \iota_X \alpha$. Indeed, Cartan's formula gives

\begin{equation}
	\label{eq:differential_K}
	\mathrm{d}K = \mathrm{d}\iota_X \alpha = \mathcal L_X \alpha - \iota_X \mathrm{d}\alpha = - \iota_X \mathrm{d}\alpha.
\end{equation}

Since $R$ is in the kernel of $\mathrm{d}\alpha$, we get that $\mathcal L_R K = \iota_R \mathrm{d}K = 0$, and hence $K$ is preserved by the flow. We give some properties of $K$: there is a relation between critical points of $K$ and orbits of the Reeb flow which are also orbits of the action. In addition, the period of those orbits is related to the corresponding critical values of $K$. Furthermore, singular orbits consist of critical points of $K$.

\begin{prop}
  	\label{prop:crit_K_S1_Reeb}
	% There is a one-to-one correspondance between critical sets of $K$ and orbits of the Reeb flow whose image on $M$ is an orbit of the action. If $x$ is a critical point of $K$ and if $x$ has stabilizer $\mathbb{Z}/ p \mathbb{Z}$ for the $\mathbb{S}^1$-action, then the Reeb orbit of $x$ has period $\frac{|K(x)|}{p}$.
	A point $x \in M$ is a critical point of $K$ if and only if the Reeb orbit of $x$ is closed and is also an orbit of the action. In this case, if $x$ has stabilizer $\mathbb{Z}/ p \mathbb{Z}$ for the $\mathbb{S}^1$-action, then the Reeb orbit of $x$ has period $\frac{|K(x)|}{p}$.
\end{prop}

\begin{proof}
	The identity (\ref{eq:differential_K}) gives that $(\mathrm{d}K)_x$ vanishes if and only $X_x$ is in the kernel of $\mathrm{d}\alpha$, which is spanned by $R$. This only happens if $X_x$ and $R_x$ are parallel. By invariance of $R$, this only happens when $X$ and $R$ are parallel along the $\mathbb{S}^1$-orbit of $x$, i.e. when the Reeb orbit of $x$ is a $\mathbb{S}^1$-orbit. If $\gamma$ is the Reeb orbit of $x$, we have
	\[
		\int_\gamma \alpha = \left| \int_{0}^{\frac{1}{p}} \alpha_{t\cdot x}(X) \mathrm{d}t\right| = \frac{|K(x)|}{p}
	\]
	which proves the result.
\end{proof}

A corollary of Proposition \ref{prop:crit_K_S1_Reeb} is that $0$ is always a regular value of $K$. Indeed, $K$ vanishes exactly when $X$ lies in the contact structure, hence $X$ and $R$ cannot be colinear when $K$ vanishes. A consequence of this is that $K^{-1}(0)$ has a finite number of connected components.

\begin{prop}
	\label{prop:singular_orbits_critical}
	Let $\gamma$ be a singular orbit of $M$, and $\alpha$ an invariant contact form. Then $\gamma$ consists of critical points of $K^\alpha$.
\end{prop}

\begin{proof}
	As explained in Section \ref{sub:seifert_bundles}, a singular orbit $\gamma$ with surgery coefficients $(p,q)$ has a local model given by $(rt_1, t_2) \in D^2 \times \mathbb{S}^1 \mapsto r t_1^p t_2^q \in D^2$. If $x \in \gamma$ is mapped to $(0,t_0)$ in this local picture, then the linearized circle action at $x$ acts at time $\frac{1}{p}$ as a rotation of angle $\frac{2 \pi}{p}$ on $T_{(0,t_0)} (D^2 \times \{t_0\})$. We write $R_{\frac{2 \pi}{p}}$ for that rotation. By $\mathbb{S}^1$-invariance of $K^\alpha$, we get that the identity $\mathrm{d}K^\alpha_{x} = \mathrm{d}K^\alpha_{x} R_{\frac{2 \pi}{p}}$ holds on $T_{(0,t_0)} (D^2 \times \{t_0\})$, meaning that $Id - R_{\frac{2 \pi}{p}}$ has image in the kernel of $(\mathrm{d}K^\alpha_{x})_{|T_{(0,t_0)} (D^2 \times \{t_0\})}$. Since $Id - R_{\frac{2 \pi}{p}}$ is invertible when $p \geq 2$, this implies that the restriction of $(\mathrm{d}K^\alpha)_{x}$ to $T_{(0,t_0)} (D^2 \times \{t_0\})$ vanishes. The circle action being tranverse to $D^2 \times \{t_0\}$, $(\mathrm{d}K^\alpha)_x$ must be zero, which proves the result.
\end{proof}

Let $\alpha$ be an invariant contact form. If the Seifert bundle $M$ has singular orbits, we write $(a_i)_{1 \leq i \leq L}$ for the singular orbits of the action. Otherwise, we write $a_1$ for an arbitrary $\mathbb{S}^1$-orbit whose image coincides with the image of a closed Reeb orbit. In both cases, we write $\mathcal A$ for the collection $\mathcal A = (a_i)_{1 \leq i \leq L}$.

We finish this section by proving the inequality stated in Theorem \ref{thm:ineg_sys} when $K$ does not vanish, i.e. when the action is transverse to the contact structure. When the Seifert bundle is a circle bundle and the action is transverse to the contact structure, this inequality was already proved in the blog post \cite{hutchings_2013} of Hutchings. The proof given in \cite{hutchings_2013} extends to the case of Seifert bundles with no difficulty.

\begin{prop}
    \label{prop:ineg_sys_facile}
	Let $\alpha \in \Omega(g,(p_i, q_i)_{1 \leq i \leq L})$ such that $K^\alpha$ does not vanish. Then the inequality
	\[
		\sys(\alpha)^2 \leq \frac{\vol(\alpha)}{|e|}
	\]
	holds. Moreover, equality holds if and only if $\alpha$ is a Zoll contact form.
\end{prop}

\begin{proof}
	Since $K^\alpha$ does not vanish, $\alpha_0 = \frac{\alpha}{K^\alpha}$ is an invariant contact form on $M$, and $K^{\alpha_0} = 1$. By Proposition \ref{prop:crit_K_S1_Reeb}, if $k$ is a critical value of $K^\alpha$, then $|k|$ is larger than $\sys(\alpha)$. In particular, the minimum of $K^\alpha$ has absolute value larger than $\sys(\alpha)$, and 
	\[
		\vol(\alpha) = \int_{M} \alpha \wedge \mathrm{d}\alpha = \int_{M} (K^\alpha)^2 \alpha_0 \wedge \mathrm{d}\alpha_0 \geq \text{min}(K^\alpha)^2 \int_{M} \alpha_0 \wedge \mathrm{d}\alpha_0 \geq \sys(\alpha)^2 \int_{M} \alpha_0 \wedge \mathrm{d}\alpha_0.
	\]
	Since $\int_{M} \alpha_0 \wedge \mathrm{d}\alpha_0 = |e|$ (see for example Proposition 6.1 in \cite{Geiges2022}), we get the claimed systolic inequality. Moreover, equality holds if and only if both inequalities are equalities. The first equality requires $K^\alpha$ to be constant, which occurs when all Reeb orbits are also $\mathbb{S}^1$-orbits. This happens exactly when $\alpha$ is Besse and the action is given by the Reeb flow of $\alpha$. In that case, if $K^\alpha$ is constant equal to $K^0$, then regular orbits are closed Reeb orbits of period $K_0$. Furthermore there is a finite number of singular orbits which have period $\frac{K_0}{p_i}$, with $p_i \in \mathbb{N}, p_i > 1$. Hence, the systole of $\alpha$ is $\frac{K_0}{\max p_i}$. As a consequence, the second inequality is an equality when there is no singular orbits, meaning that $\alpha$ is Zoll.
\end{proof}

% subsection invariant_contact_forms_and_first_properties (end)

\subsection{Classification of invariant contact structures} % (fold)
\label{sub:classification_of_invariant_contact_structures}

We recall that (coorientable) contact structures are globally defined as the kernel of a contact form. A contact structure is \textbf{invariant} if it is induced by an invariant contact form, or equivalently if it is invariant under the action as a plane field. In this section we provide without proofs the classification of invariant contact structures on Seifert bundles. For a more detailed account of the classification, see \cite{Lutz1977} and \cite{Giroux1999StructuresDC} in the case of a free $\mathbb{S}^1$-action, \cite{niederkruger2005} for the general case of Seifert bundles, as well as \cite{GeigesHedickeSaglam} for a detailed discussion of $\mathbb{S}^1$-invariant contact structures on $\mathbb{S}^3$. Invariant contact structures are classified by their set of Legendrian $\mathbb{S}^1$-orbits. The following definition associates to an invariant contact structure a graph, encoding the relevant properties of those Legendrian orbits.

\begin{definition}
	Let $\xi$ be an invariant contact structure, defined as the kernel of an invariant contact form $\alpha$. We define a graph $G(\alpha) = (V,E)$ as follows: vertices of $G(\alpha)$ are connected components of $M \setminus (K^\alpha)^{-1}(0)$, and there is an edge between two vertices when the corresponding components of $M \setminus (K^\alpha)^{-1}(0)$ share a boundary component. If a component of $M \setminus (K^\alpha)^{-1}(0)$ contains singular orbits, we label the corresponding vertex with the surgery coefficients $(p_i,q_i)$ of those singular orbits. Since this labelled graph does not depend on the choice of the invariant contact form $\alpha$ inducing $\xi$, we denote it by $G(\xi)$.
\end{definition}

The classification of invariant contact structures can now be stated. First, a uniqueness result states that the labelled graph $G(\xi)$ uniquely determines, up to equivariant diffeomorphism, $\xi$ as an $\mathbb{S}^1$-invariant contact structure.

\begin{prop}[Section 2 in \cite{Lutz1977}, Lemma IV.10 in \cite{niederkruger2005}]
	Let $\xi_1, \xi_2$ be two postive invariant contact structures on $M$. There is an equivariant diffeomorphism $\phi$ of $M$ such that $\phi_* \xi_1 = \xi_2$ if and only if $G(\xi_1)$ and $G(\xi_2)$ are isomorphic.
\end{prop}

Second, an existence result essentially describes the class of graphs in the range of the map $\xi \mapsto G(\xi)$. Since $0$ is a regular value of $K^\alpha$, $K^\alpha$ must have different signs on any two components of $M \setminus (K^\alpha)^{-1}(0)$ which share a boundary component. As a result, the graph $G(\xi)$ is always bipartite, where a graph $G$ is \textbf{bipartite} when there is a labelling of the vertices of $G$ over a set with two elements such that any two vertices connected by an edge have different labels. Given $S$ a closed oriented surface, we say that a bipartite graph $G$ is \textbf{S-admissible} if there is a non-empty collection $\mathcal C$ of embedded circles in $S$ and a bijection between vertices of $G$ and components of $S \setminus \mathcal C$ such that there is an edge between two vertices if and only if the corresponding components of $S \setminus \mathcal C$ share a boundary component. We insist on the fact that the graph with one vertex and no edges is not considered as S-admissible, and is treated separately.

\begin{prop}[Section 1.4 in \cite{Lutz1977}, Lemma IV.15 in \cite{niederkruger2005}]
	Let $S$ be a closed oriented surface with genus $g$ and $G$ a S-admissible graph. Assume that the vertices of $G$ are labelled by a (possibly empty) finite set of coprime integers $(p_i,q_i)$. Then there exists a positive $\mathbb{S}^1$-invariant contact structure $\xi$ on the Seifert bundle $M(g, (p_i, q_i)_{1 \leq i \leq L})$ whose associated labelled graph $G(\xi)$ is $G$.
\end{prop}

The following result deals with the case of the graph with one vertex and no edges. This graph corresponds to invariant contact structures which are transverse to the orbits of the $\mathbb{S}^1$-action. The study of contact structures transverse to a circle action was started in \cite{bb725b3c-2f5a-3ef1-b918-e4d2290d3cc4}. An invariant contact structure on a Seifert bundle $M$ is called \textbf{positive} if the orientation it induces on $M$ coincides with the orientation of $M$ as a Seifert bundle.

\begin{prop}[Section 1.3 in \cite{Lutz1977}, Theorem 3 in \cite{bb725b3c-2f5a-3ef1-b918-e4d2290d3cc4}, Proposition 3.1 in \cite{Lisca2004}]
	\label{prop:existence_Zoll}
	If a Seifert bundle $M$ satisfies $e < 0$, then there is a positive contact structure on $M$ which is transverse to the fibers of the bundle.
\end{prop}
In that case, the Seifert fibration can be realized as the Reeb flow of a Besse contact form, see Theorem 1.4 in \cite{Kegel2021}.

\subsection{Examples} % (fold)
\label{sub:examples}

We now give examples of invariant contact forms.

\begin{enumerate}
	\item \label{example:zoll} \textbf{Besse and Zoll contact forms, Boothby-Wang construction: }
	A contact form $\alpha$ is called \textbf{Besse} if all its Reeb orbits are closed. If in addition all Reeb orbits have the same minimal period, $\alpha$ is said to be \textbf{Zoll}. The Reeb flow of a Besse (resp. Zoll) contact form induces an almost-free (resp. free) $\mathbb{S}^1$-action by strict contactomorphisms. This action is transverse to the contact structure, in other words the graph $G(\alpha)$ has one vertex and no edges. The existence of Zoll contact forms is described by a classical result of Boothby--Wang \cite{bb725b3c-2f5a-3ef1-b918-e4d2290d3cc4}: they are in one-to-one correspondence with principal circle bundles over a base manifold $B$, with curvature given by an integral symplectic form on $B$. A similar statement has been proved in \cite{Kegel2021} for Besse contact forms, where $B$ is now an orbifold. Finally, isotopy classes of Zoll contact forms have been fully described in dimension three in \cite{benedetti_kang_2020}.

	\item \label{example:metrics_revolution} \textbf{Metrics of revolution on $\mathbb{S}^2$: }A natural source of examples for contact forms arise from metric geometry. It is a classical fact that a Riemannian metric, or more generally a Finsler metric, on a manifold induces a canonical contact form on its unit tangent bundle. In our case, we can produce $\mathbb{S}^1$-invariant contact forms from metrics as follows: any effective circle action on $\mathbb{S}^2$ has exactly two fixed points, and lifts to a free circle action on the unit tangent bundle of $\mathbb{S}^2$, which is diffeomorphic to $\mathbb{RP}^3$. This action makes $\mathbb{RP}^3$ the total space of a principal bundle over $\mathbb{S}^2$ with Euler number $2$. Any Finsler metric $F$ on $\mathbb{S}^2$ which is invariant under the action gives rise to a contact form $\alpha_F$ on $\mathbb{RP}^3$ which is invariant under the lifted free action. The set of Legendrian orbits for the lifted action consists of the union of the geodesics on $\mathbb{S}^2$ joining the two fixed points of the action. In particular, this set is connected, and the graph $G(\alpha_F)$ consists of two vertices connected by a single edge. Similarly, the metrics of revolution on spindle orbifolds considered in \cite{MR4529024} induce $\mathbb{S}^1$-invariant contact forms on lens spaces other than $\mathbb{RP}^3$.

	\item \textbf{Multiplication by an invariant function: }Multiplying an invariant contact form by an invariant nowhere-vanishing function produces a new invariant contact form. This can of course change drastically the Reeb flow, and does not affect the graph $G(\alpha)$.

	\item \textbf{Diverse surgeries:}
	\begin{itemize} 

		\item \textbf{Lutz twists: }A useful topological operation on contact structures is the Lutz twist, a topologically trivial surgery producing a contact structure in a different homotopy class from the original one. We describe that operation in the invariant setting (see \cite{geiges_2008} for more details). Let $\alpha$ be an invariant contact form. Since $K^\alpha$ always has a critical point, $\alpha$ has a closed Reeb orbit $\gamma_0$ which coincides with an $\mathbb{S}^1$-orbit. Up to multiplying $\alpha$ with an appropriate invariant function, we can assume that $K^\alpha$ is constant in a neighborhood of $\gamma_0$, meaning that $\gamma_0$ has a neighborhood $U$ diffeomorphic to a solid torus $\mathbb{S}^1\times \mathbb{D}$, all of whose $\mathbb{S}^1$-fibers are both $\mathbb{S}^1$-orbits and closed Reeb orbits, and oriented in the Reeb direction. In suitable coordinates, the restriction of $\alpha$ to $U$ has the form $\alpha = K_0\mathrm{d}\theta + r^2 \mathrm{d}\phi$, with $\theta$ the $\mathbb{S}^1$ coordinate on $\mathbb{S}^1\times \mathbb{D}$, and $(r, \phi)$ polar coordinates on $\mathbb{D}$. We fix two functions $h_1, h_2 \colon [0,1] \rightarrow \mathbb{R}$ satisfying the following properties:
	    \begin{enumerate}
	    	\item $h_1(r) = -K_0$ and $h_2(r) = -r^2$ for $r$ close to $0$,
	    	\item $h_1(r) = K_0$ and $h_2(r) = r^2$ for $r$ close to $1$,
	    	\item $(h_1(r), h_2(r))$ is never parallel to $(h_1'(r), h_2'(r))$,
	    \end{enumerate}
	    and define $\tilde \alpha = h_1(r) \mathrm{d}\theta + h_2(r) \mathrm{d}\phi$ on $U$. Since $\tilde \alpha$ coincides with $\alpha$ in a neighborhood of $\partial U$, replacing $\alpha$ with $\tilde \alpha$ on $U$ defines a new contact form, the contact condition being ensured by the last condition on $(h_1, h_2)$. Moreover, since $h_1$ and $h_2$ are $\mathbb{S}^1$-invariant on $U$, $\tilde \alpha$ is $\mathbb{S}^1$-invariant on $M$. We have that $K^{\tilde \alpha}$ coincides with $K^\alpha$ outside of $U$, and is given by $h_1$ on $U$. While $K^\alpha$ is constant on $U$, $h_1$ does vanish on $U$, with vanishing set a two-torus cutting $U$ into a smaller solid torus and its complement. As a consequence, the graph $G(\tilde \alpha)$ is $G(\alpha)$ with a new vertex and an edge between this vertex and the component where the Lutz twist is performed. Iterating this operation any finite number of times allows to construct from $\alpha$ invariant contact forms whose associated graph is $G(\alpha)$ concatenated with a finite number of finite trees.

		\item \textbf{Connected sums: } Generalizing the previous construction, one can also perform connected sums of invariant contact structures. For $i \in \{1,2\}$ let $\alpha_i$ an invariant contact form on $M_i$, admitting a solid torus $U_i \subset M_i$ fibered by closed Reeb orbits of period 1, all coinciding with orbits of the action. As previously, this can be achieved by multiplying by an appropriate invariant function. In suitable coordinates, the restriction of $\alpha_i$ to $U_i$ has the form $\alpha_i = \mathrm{d}\theta + r^2 \mathrm{d}\phi$, with $\theta$ the $\mathbb{S}^1$ coordinate on $U_i \simeq \mathbb{S}^1\times \mathbb{D}$, and $(r, \phi)$ polar coordinates on $\mathbb{D}$. Gluing $M_1$ to $M_2$ by identifying $U_1$ and $U_2$ via those coordinates produces a free $\mathbb{S}^1$-action on the connected sum, as well as a global contact form $\alpha$ invariant under that action. The graph $G^\alpha$ associated to $\alpha$ is the graph obtained after identifying the vertices of $G^{\alpha_1}$ and $G^{\alpha_2}$ corresponding to the components of $M_1 \setminus (K^{\alpha_1})^{-1}(0)$ and $M_2 \setminus (K^{\alpha_2})^{-1}(0)$ on which the sum was performed. Iterating this construction with different choices of $M_i$ and $\alpha_i$ allows to show that $\Omega(S,e,G)$ is non-empty when $G$ is admissible.

		\item \textbf{Introducing singular orbits: } Finally, we can exit the realm of principal circle bundles and enter that of Seifert fibrations by introducing singular orbits. To do so, let $(p,q)$ be two coprime integers, and $\alpha \in \Omega(g,0)$ an invariant contact form on $\mathbb{S}^1 \times S_g$. We can assume as before that $K^\alpha$ is constant on an embedded solid torus $U$. There is a free $\mathbb{Z}/ p \mathbb{Z}$ action on $U \simeq \mathbb{S}^1 \times \mathbb{D}$, generated by $1 \cdot (t,z) = \left(t+\frac{q}{p}, e^{\frac{2i \pi}{p}}z\right)$, which leaves $\alpha$ invariant. The quotient $U/\left(\mathbb{Z} / p \mathbb{Z}\right)$ is diffeomorphic to a solid torus, carrying an almost free circle action which has a unique singular orbit with isotropy $\mathbb{Z}/ p \mathbb{Z}$. Furthermore, $\alpha$ induces an invariant contact form on $U/\left(\mathbb{Z} / p \mathbb{Z}\right)$. After gluing back $U/\left(\mathbb{Z} / p \mathbb{Z}\right)$ to $M \setminus U$ as described in Section \ref{sub:seifert_bundles}, we obtain an invariant contact form on $M(g,(p,q))$. The effect of this operation on $G^\alpha$ is to add a label $(p,q)$ on the vertex corresponding to the component of $M \setminus (K^\alpha)^{-1}(0)$ containing $U$. Iterating the process allows to introduce any finite number of singular orbits.
	\end{itemize}

\end{enumerate}

% subsection examples (end)

\section{Invariant surfaces of section and potentials}
\label{sec:surface_of_section_and_potentials}
This section contains the two main ingredients for the proof of Theorem \ref{thm:ineg_sys}, namely an invariant surface of section, and a decomposition for invariant contact forms. This will enable us to introduce a family of potentials encoding most of the dynamic of the Reeb flow. We will end this section by giving some important properties of those potentials.

\subsection{An invariant surface of section}
\label{sub:invariant_surface_of_section}

Let us fix $\alpha$ in $\Omega(g, (p_i, q_i)_{1 \leq i \leq L})$. Being invariant, $K$ induces a well-defined function $\tilde K \colon M / \mathbb{S}^1 \rightarrow \mathbb{R}$. Let $\mathcal C$ be the union of the connected components of levels sets of $K \colon M \rightarrow \mathbb{R}$ that contain a critical point of $K$. By Proposition \ref{prop:singular_orbits_critical}, singular orbits of the action consist of critical points of $K$. Hence, $\mathbb{S}^1$ acts freely on $M \setminus \mathcal C$, and the quotient $\Gamma = (M \setminus \mathcal C)/\mathbb{S}^1$ is a union of cylinders $\Gamma = \bigcup \limits_{i \in I} \Gamma_i$. We will denote by $I_0$ the set of indices $i \in I$ such that $\tilde K$ vanishes on $\Gamma_i$. Since $K^{-1}(0)$ is finite, $I_0$ is finite as well. We write $\epsilon_i^+$ and $\epsilon_i^-$ for the two components of $\partial \Gamma_i$, with $\tilde K(\epsilon_i^+) > \tilde K(\epsilon_i^-)$, and $K_i^+$ (resp. $K_i^-$) for $\tilde K(\epsilon_i^+)$ (resp. $\tilde K(\epsilon_i^-)$). We choose $x_i^+$ (resp. $x_i^-$) an arbitrary point of $\epsilon_i^+$ (resp. $\epsilon_i^-$). Finally, we choose a curve $\tilde \gamma_i \colon (0,1) \rightarrow \Gamma_i$ which is transverse to the level sets of $\tilde K$, converging to $x_i^-$ at $0$ and to $x_i^+$ at $1$.

Writing $\gamma_i$ for a lift of $\tilde \gamma_i$ to $M$, we set $\Sigma_i$ as the image of $\gamma_i$ under the $\mathbb{S}^1$-action, i.e. $\Sigma_i = \mathbb{S}^1 \cdot \gamma_i = \pi^{-1}(\tilde \gamma_i)$. The $\Sigma_i$ are pairwise disjoint cylinders, whose boundary components have the form $\mathbb{S}^1 \cdot x$, with $x$ a critical point of $K$, and as such are both closed Reeb orbits and $\mathbb{S}^1$-orbits. We write $\Sigma$ for their union $\bigcup \limits_{i \in I} \Sigma_i$. The following proposition states that $\Sigma$ is a surface of section.

\begin{prop}
	\label{prop:temps_de_retour_fini}
    The vector field $R$ is transverse to $\Sigma$. Moreover, for all $x$ in $\Sigma$, the Reeb orbit of $x$ intersects $\Sigma$ in positive and negative time.
\end{prop}

\begin{proof}
    Let $x \in \Sigma$. Since the vector field $R$ is nowhere parallel to $X$ on $M \setminus \mathcal C$, its pushforward $\pi_{*}R$ is a nowhere vanishing vector field on $\Gamma$. Moreover $\pi_{*}R$ is tangent to the level sets of $\tilde K$, as $K$ is preserved by the Reeb flow. Since $\tilde \gamma_i$ is transverse to the level sets of $\tilde K$, it is transverse to $\pi_*(R)$. Thus $R$ is transverse to $\Sigma_i = \pi^{-1}(\gamma_i)$. The level sets of $\tilde K$ in $\Gamma$ are unions of circles, hence the flow lines of $\pi_* R$ are circles as well. In particular, they pass infinitely many times through $\pi(x)$, which means that the Reeb orbit of $x$ in $M$ must intersect the $\mathbb{S}^1$-orbit of $x$ infinitely many times both in negative and positive time. This proves the result. 
\end{proof}

For $x$ in $\Sigma$, we will write $\tau(x)$ for the first positive time at which the Reeb orbit starting at $x$ hits $\Sigma$ again. This allows us to define a diffeomorphism $\Phi$ of $\Sigma$ as
\[
	\begin{array}{ccccc}
		\Phi & \colon & \Sigma & \rightarrow & \Sigma \\
		& & x & \mapsto & \phi^{\tau(x)}(x)
	\end{array}
\]
as well as $\Phi_i = \Phi_{|\Sigma_i}$. Note that periodic points of $\Phi$ correspond to closed Reeb orbits of $\alpha$. We will write $c_x$ for the portion of Reeb orbit $\phi^{[0,\tau(x)]}(x)$ starting at $x$ and ending at $\Phi(x)$. If $k \in (K^-_i, K^+_i)$, there is a unique $x_{i,k}$ in $\gamma_i$ such that $K(x_{i,k}) = k$. We will occasionally write $c_{i,k}$ for $c_{x_{i,k}}$.

\subsection{A decomposition for invariant contact forms} % (fold)
\label{sub:a_decomposition_for_invariant_contact_forms}

In this section, we give a decomposition for invariant contact forms on $M$. We start with a preliminary lemma.

\begin{lemma}
    \label{lemma:lift_form_quotient}
    Let $U$ be a $\mathbb{S}^1$-invariant open set of $M$ on which $\mathbb{S}^1$ acts freely, and $\kappa$ a one-form on $U$. Then there is a one-form $\beta$ on $\pi(U)$ such that $\kappa = \pi^* \beta$ if and only if $\mathcal L_{X} \kappa = 0$ and $i_{X} \kappa = 0$.
\end{lemma}

\begin{proof}
	The direct implication is easy. Conversely, let $\kappa$ be such that $\mathcal L_{X} \kappa = 0$ and $i_{X}\kappa = 0$. We define $\beta$ by $\beta_x(\tilde Y) = \kappa_y(Y)$, where $\pi(y)=x$ and $\mathrm{d} \pi_y(Y) = \tilde Y$. Since $\mathcal L_{X}\kappa = 0$, $\beta_y(Y)$ does not depend on the choice of $y$ in $\pi^{-1}(x)$. Neither does it depend on the choice of $Y$ since $i_{X} \kappa = 0$.
\end{proof}

Recall that the invariants $(p_i,q_i)_{1 \leq i \leq L}$ are assumed to be normalized and that $\mathcal A$ is the collection of all the singular orbits of the circle action if there are any, and of an arbitrary closed Reeb orbit in the critical set of $K$ otherwise. The restriction of the circle action to $M \setminus \mathcal A$ is free, and makes it the total space of a trivial circle bundle over a punctured surface $\tilde S = S_g \setminus \{u_1, \ldots, u_L\}$. Here, each puncture $u_i$ corresponds to the boundary component $a_i$ of $M \setminus \mathcal A$. It will be useful to choose a section $S$ of this bundle as follows. We explained in Section \ref{sub:seifert_bundles} how $M(g, (p_i,q_i)_{1 \leq i \leq L})$ is obtained by performing a finite number of Dehn surgeries with coefficients $(p_i,q_i)$ on $\mathbb{S}^1 \times S_g$ along the unknots $\mathbb{S}^1 \times \{u_i\}$. Fixing $t_0 \in \mathbb{S}^1$, then the image $S$ in $M$ of $\{t_0\} \times \tilde S$ after surgery is a section for $M \setminus \mathcal A \rightarrow \tilde S$. This section has the following property, which will be useful in the proof of Lemma \ref{lemma:computing_euler}: if $b_i$ is the boundary component of $S$ in $M$ corresponding to the puncture $u_i$, given the boundary orientation, then $b_i$ covers $-q_i$ times $a_i$.

This allows to define a closed one-form $\sigma$ on $M \setminus \mathcal A$ such that:
\begin{itemize}
    \item $\sigma$ vanishes on $\pi_1(S)$ ;
    \item $\int_{\mathbb{S}^1\cdot x} \sigma = 1$ for all $x$ in $M \setminus \mathcal A$ ;
    \item $\sigma$ is $\mathbb{S}^1$-invariant.
\end{itemize}
Now, the one-form $\alpha - K \sigma$ defined on $M \setminus \mathcal A$ is $\mathbb{S}^1$-invariant, and satisfies
\[
	\iota_X(\alpha - K \sigma) = K - K = 0.
\]
Hence, Lemma \ref{lemma:lift_form_quotient} allows to write on $M \setminus \mathcal A$
$$
    \alpha = K \sigma + \pi^*\beta,
$$
where $\beta$ is a one-form on $\tilde S$. In the next lemma, we express the contact volume $\alpha \wedge \mathrm{d}\alpha$ on $M \setminus \mathcal A$ in term of $K$, $\sigma$ and $\beta$.

\begin{lemma}
	\label{lemma:volume_form_decomposition}
	On $M \setminus \mathcal A$, we have
	\[
		\alpha \wedge \mathrm{d}\alpha = \sigma \wedge \pi^* \left(\tilde K \mathrm{d}\beta + \beta \wedge \mathrm{d}\tilde K \right).
	\]
\end{lemma}

\begin{proof}
	The decomposition $\alpha = K \sigma + \pi^* \beta$ on $M \setminus \mathcal A$ implies
	\[
		\alpha \wedge \mathrm{d}\alpha = (K \sigma + \pi^* \beta) \wedge \left(\mathrm{d}K \wedge \sigma + \pi^* \mathrm{d}\beta\right).
	\]
	Since $\sigma \wedge \sigma$ and $\pi^*\beta \wedge \pi^* \mathrm{d}\beta$ both vanish, we get that
	\[
		\alpha \wedge \mathrm{d}\alpha = \sigma \wedge \pi^* \left(\tilde K \mathrm{d}\beta + \beta \wedge \mathrm{d}\tilde K \right).
	\]
\end{proof}

% subsection a_normal_form_for_invariant_contact_forms (end)

\subsection{A family of potentials} % (fold)
\label{sub:a_family_of_potentials}
For each $i \in I$ we define the map
\[
	\begin{array}{ccccc}
		J_i & \colon & (K_i^-, K_i^+) & \rightarrow & \mathbb{R} \\
		& & k & \mapsto & \int_{\pi(c_{i,k})} \beta
	\end{array}
\]
which is well-defined as $M \setminus \mathcal C$ is contained in $M \setminus \mathcal A$.

\begin{lemma}
    \label{lemma:derivative_J}
    The map $J_i$ is differentiable on $(K_i^-, K_i^+)$, with $J_i'(k) = -\int_{c_{i,k}} \sigma$.
\end{lemma}

\begin{proof}
    Let $k^- < k^+$ be in $(K_i^-, K_i^+)$. Define $H$ to be $\bigcup \limits_{k^- < k < k^+}c_{i,k}$. This surface projects injectively on $\Gamma$, with image $\Gamma_i \cap K^{-1}((k^-,k^+)$, whose oriented boundary is $\pi(k^-) - \pi(k^+)$ (both oriented by the flow of $\pi_* R$). By Stokes theorem, one then has
    \begin{align*}
        J_i(k^+)-J_i(k^-) &= -\int_{\pi(H)} \mathrm{d}\beta \\
        &= -\int_{H} \pi^*\mathrm{d}\beta \\
        &= -\int_{H} \left(\mathrm{d}\alpha - \mathrm{d}K \wedge \sigma\right) \\
        &= -\int_{H} \sigma \wedge \mathrm{d}K \\
        &= \int_{k^-}^{k^+}\left(-\int_{c_{i,k}} \sigma \right) \mathrm{d}k.
    \end{align*}
    The penultimate equality follows from the fact that the restriction of $\mathrm{d} \alpha$ to $TH$ vanishes as $R$ is tangent to $H$, and the last equality is obtained after the change of variable $x \in \gamma_i \cap H \mapsto K(x) \in (k^-, k^+)$. Dividing by $k^+ - k^-$ and letting $k^+$ go to $k^-$ yields the result.
\end{proof}

An immediate corollary of Lemma \ref{lemma:derivative_J} is the identity below relating $\tau$ and $J_i$.

\begin{coro}
    \label{coro:identity_tau}
    For all $k \in (K_i^-, K_i^+)$, $\tau(x_{i,k}) = J_i(k) - k J_i'(k)$. In particular, $J_i(k) - k J_i'(k) > 0$.
\end{coro}

\begin{proof}
    Integrating the identity $\alpha = K \sigma + \pi ^* \beta$ along $c_{i,k}$ yields
    $$
        \tau(x_{i,k}) = k \int_{c_{i,k}}\sigma + J_i(k).
    $$
    By Lemma \ref{lemma:derivative_J}, $\int_{c_{i,k}}\sigma = - J_i'(k)$ and the result follows.
\end{proof}

A second corollary of Lemma \ref{lemma:derivative_J} is that the derivative of $J_i$ is closely related to periodic points of the first-return map of the surface of section.

\begin{coro}
    \label{coro:periodic_point}
    The Reeb orbit starting at $x_{i,k} \in \gamma_i$ is closed if and only if $J_i'(k)$ is rational. Moreover, if $J_i'(k) = \frac{p}{q}$ with $\gcd(p,q)=1$, then $x_{i,k}$ is a $q$-periodic point of $\Phi_i$, and the corresponding closed Reeb orbit has minimal period $q(J_i(k) - k J_i'(k))$.
\end{coro}

\begin{proof}
    Since $K$ is preserved by the Reeb flow, the first-return map $\Phi_i$ of $\Sigma_i$ induces a circle diffeomorphism of $\mathbb{S}^1 \cdot x_{i,k}$ for each $k \in (K_i^-, K_i^+)$. By $\mathbb{S}^1$-invariance of the Reeb flow, each of those circle diffeomorphisms is a translation. Finally, homotoping $c_{i,k}$ relative to its boundary in $M \setminus \mathcal A$ to a path contained in $\mathbb{S}^1 \cdot x_{i,k}$, together with the closedness of $\sigma$, gives us that the circle diffeomorphism of $\mathbb{S}^1 \cdot x_{i,k}$ is a translation of shift $\int_{c_{i,k}} \sigma$. As a consequence, the Reeb orbit starting at $x_{i,k}$ is closed if and only if $\int_{c_{i,k}} \sigma = -J_i(k)$ is a rational number.

    When $J_i'(k) = \frac{p}{q}$ with $\gcd(p,q)=1$, the translation of shift $\int_{c_{i,k}} \sigma = - \frac{p}{q}$ is $q$-periodic. Furthermore, the closed Reeb orbit starting at $x_{i,k}$ has period $q \tau(x_{i,k}) = q(J_i(k) - k J_i'(k))$ according to Corollary \ref{coro:identity_tau}.
\end{proof}

Finally, the next lemma states how to compute the contact volume from the potentials.

\begin{lemma}
	\label{lemma:volume_J}
	The contact volume satisfies
	\[
		\vol(M) \geq \sum_{i \in I} \int_{K_i^-}^{K_i^+}\left(J_i(k) - k J_i'(k)\right) \mathrm{d}k
	\]
	and each term in this sum is positive.
\end{lemma}

\begin{proof}
	Let $U = \bigcup \limits_{x \in \Sigma} \phi^{[0, \tau(x)]}(x)$ be the image of the surface of section $\Sigma$ by the Reeb flow. Since the Reeb segments $\phi^{[0, \tau(x)]}(x)$ are tangent to the kernel of $\mathrm{d}\alpha$, we have
	\[
		\vol(M) = \int_M \alpha \wedge \mathrm{d}\alpha \geq \int_{U} \alpha \wedge \mathrm{d}\alpha = \int_{\Sigma}\tau \mathrm{d}\alpha.
	\]
	Furthermore, $\mathrm{d}\alpha = \mathrm{d}K \wedge \sigma + \pi^*\beta$ on $U$. Using that $\pi^* \beta$ vanishes on $\Sigma$, and $\Sigma = \mathbb{S}^1 \cdot \bigcup_{i \in I} \gamma_i$, we get that
	\[
		\int_{\Sigma}\tau \mathrm{d}\alpha = \sum \limits_{i \in I} \int_{\gamma_i}\tau \mathrm{d}K = \sum \limits_{i \in I} \int_{K_i^-}^{K_i^+}\tau(x_{i,k}) \mathrm{d}k
	\]
	where the last equality follows from the change of variable $x \in \gamma_i \mapsto K(x) \in (K_i^-, K_i^+)$. Finally, Corollary \ref{coro:identity_tau} yields the claimed inequality.
\end{proof}

We finish this section by giving a relation between the Euler number of the Seifert bundle and the one-form $\beta$.

\begin{lemma}
	\label{lemma:computing_euler}
	Let $\alpha \in \Omega\left(g,(p_l,q_l)_{1 \leq l \leq L}\right)$. Then
	\[
		\int_{\partial \tilde S} \frac{\beta}{\tilde K} = e.
	\]
\end{lemma}

\begin{proof}
	For $1 \leq l \leq L$, the boundary component $b_l$ of $S$ is the lift of a boundary component $\tilde b_l$ of $\tilde S$. Recall that $S$ was chosen so that $b_l$ covers $-q_l$ times $a_l$, $b_l$ being oriented as a boundary and $a_l$ as an orbit of the $\mathbb{S}^1$-action. Let $K_l$ be the value of $K$ on $a_l$. We fix $(b_{l,n})_{n \in \mathbb{N}}$ a sequence of compact one-dimensional submanifolds in the interior of $S$, converging to $b_l$ in $H_1(S,\mathbb{Z})$.

	By definition of $\sigma$ as a connection one-form, $\sigma$ vanishes on $\pi_1(\tilde S)$, in particular on the curves $b_{l,n}$. It follows that 
	\[
		\int_{b_{l,n}} \frac{\pi^*\beta}{K} = \int_{b_{l,n}} \frac{\alpha}{K}.
	\]
	Since $b_l = -q_l a_l$ in $M$, we get that $\int_{b_{l,n}} \frac{\alpha}{K}$ converges to $-q_l \frac{\int_{a_l} \alpha}{K_l}$. If $K_l > 0$, then $a_l$ is positively oriented as a Reeb orbit. Hence,
	\[
		-q_l \frac{\int_{a_l} \alpha}{K_l} = -\frac{q_l K_l}{p_l K_l} = - \frac{q_l}{p_l}.
	\]
	Similarly, $a_l$ is negatively oriented as a Reeb orbit when $K_l < 0$, and we have
	\[
		-q_l \frac{\int_{a_l} \alpha}{K_l} = \frac{q_l |K_l|}{p_l K_l} = - \frac{q_l}{p_l}.
	\]
	The result follows when summing over $1 \leq l \leq L$.
\end{proof}

\section{Necessary properties of contact forms with high systolic ratio}
\label{sec:necessary_properties_high_systolic_ratio}

The goal of this section is to prove Proposition \ref{prop:potential_high_systolic_ratio_is_almost_linear} below, which states that if an hypothetical contact form $\alpha \in \Omega\left(g,(p_l,q_l)_{1 \leq l \leq L}\right)$ has high systolic ratio, then its potentials $(J_i)_{i \in I_0}$ must be uniformly close to linear maps. Throughout this section, $J \colon [-1,1] \rightarrow \mathbb R$ will be a smooth function satisfying the following conditions:

\begin{enumerate}[label=$(C_{\arabic*})$]
    \item $\forall x \in (-1,1), J(x) - x J'(x) > 0$;
    \label{enumitem:cond_high_systolic_ratio_contact}
    \item If $x$ is such that $J'(x)$ is rational, with $J'(x) = \frac{p}{q}$, and $q > 0$, then $J(x) - x J'(x) \geq \frac{1}{q}$;
    \label{enumitem:cond_high_systolic_ratio_systole}
    \item $\int_{-1}^1 (J(t)-tJ'(t)) \mathrm{d}t < \frac{1}{80}$.
    \label{enumitem:cond_high_systolic_ratio_volume}
\end{enumerate}

We start with two preliminary lemmas which give respective consequences of \ref{enumitem:cond_high_systolic_ratio_contact} and \ref{enumitem:cond_high_systolic_ratio_systole}.

\begin{lemma}
	\label{lemma:consequence_volume}
	Let $f \colon [0,1] \rightarrow [0, +\infty[$ be a smooth function satisfying 
	\begin{equation}
		\label{enumitem:cond_f_3}
		\forall x \in (0,1), f(x) - x f'(x) > 0.
	\end{equation}
	Then the upper bound
	\[
		f(x) < \frac{2\int_0^1f(t) \mathrm{d}t}{x}
	\]
    holds.
\end{lemma}

\begin{proof}
	Since $\left(\frac{f(x)}{x} \right)'= \frac{f'(x)x - f(x)}{x^2}$, condition (\ref{enumitem:cond_f_3}) implies that the function $x \mapsto \frac{f(x)}{x}$ is decreasing. Hence, for all $0 < t < x$, one has $t \frac{f(x)}{x} < f(t)$. Integrating this inequality with respect to $t$ over $[0,x]$ yields that $x f(x) < 2 \int_0^{x}f(t) \mathrm{d}t$. Using also $f \geq 0$, one has
	\[
	 	x f(x) < 2 \int_0^1f(t) \mathrm{d}t
	\]
	for all $x$ in $[0,1]$.
\end{proof}

\begin{lemma}
    \label{lemma:consequence_systole}    
    Let $f \colon \left[0, \frac{1}{2}\right] \rightarrow \mathbb R$ be a smooth function such that if $x$ is such that $f'(x)$ is rational, with $f'(x) = \frac{p}{q}$, and $q > 0$, then
    \begin{equation}
    	\label{enumitem:cond_f_4}
    	f(x) - x f'(x) \geq \frac{1}{q}
    \end{equation}
	Assume furthermore that there exists $a \in \left(0, \frac{1}{20}\right)$ and $b \in \mathbb{R}$ such that $\max \limits_{x \in \left[\frac{1}{4}, \frac{1}{2}\right]}|f(x) - b x| < a$. Then we have that
	\[
		\max \limits_{x \in \left[0, \frac{1}{4}\right]} |f'(x) - b| < 28a.
	\]
\end{lemma}

\begin{proof}
	The function $g \colon \left[0, \frac{1}{2}\right] \rightarrow \mathbb{R}$ defined by $g(x) = f(x)-bx$ satisfies
	\[
		\max \limits_{x \in \left[\frac{1}{4}, \frac{1}{2}\right]} |g(x)| < a.
	\]
	In particular, both $\left|g(\frac{1}{4})\right|$ and $\left|g(\frac{1}{2})\right|$ are less than $a$. By the mean value theorem, there exists $x \in \left[\frac{1}{4}, \frac{1}{2}\right]$ such that $g'(x) = 4\left(g\left(\frac{1}{2}\right) - g(\frac{1}{4})\right)$. It follows that the set $\{x \in \left[\frac{1}{4}, \frac{1}{2}\right], |g'(x)|\leq 8a\}$ is not empty, and we fix $x_0$ an element of this set.

	Let us proceed by contradiction and assume that there is $x_1 \in \left[0, \frac{1}{4}\right]$ such that $|g'(x_1)| \geq 28a$. This implies that the interval $[8a, 28a]$ is contained in the image of $|g'|$, and thus that the image of $[x_1, x_0]$ by $|f'|$ contains an interval $I$ of length $20a$.

	By the pigeonhole principle, any interval of length $l < 1$ contains a rational number $\frac{p}{q}$ with $\frac{1}{q} \geq \frac{l}{2}$. Since $a$ is less than $\frac{1}{20}$, it follows that $I$ contains a rational number $\frac{p}{q_0}$ with $\frac{1}{q_0} \geq 10a$. Let $p_0 = \lfloor q_0 f'(x_0) \rfloor$, and $p_1 = p_0 + 1$. Since the image of $f'$ on $[x_1,x_0]$ contains both $f'(x_0)$ and a number of the form $\frac{p}{q_0}$, and since $f'$ is continuous, the set
	\[
		\left\{x \in \left[x_1, x_0\right], f'(x) = \frac{p}{q_0}, p \in \{p_0,p_1\}\right\}
	\]
	is not empty. Let $x_2$ be its maximum. We then have that the image of $(x_2,x_0)$ by $f'$ is an interval of length at most $\frac{1}{q_0}$. For all $x \in (x_2, x_0)$, we have
	\begin{align}
		|g'(x)| &\leq |g'(x_0)| + |g'(x_0) - g'(x)| \nonumber \\
		&=|g'(x_0)| + |f'(x_0) - f'(x)| \nonumber \\
		\label{align:borne_g_prime}
		&< 8a + \frac{1}{q_0}.
	\end{align}
	This implies that
	\begin{equation}
		\label{align:borne_g}
		|g(x_0) - g(x_2)| = \left|\int_{x_2}^{x_0} g'(t)\mathrm{d}t\right| < (x_0 - x_2)\left(8a + \frac{1}{q_0}\right).
	\end{equation}

	From (\ref{align:borne_g}), we deduce that
	\begin{equation}
		\label{eq:borne_f_x_2}
		\left| f\left(x_2\right) - bx_2 \right| = |g(x_2)| \leq |g(x_0)| + |g(x_0) - g(x_2)| < a + (x_0-x_2) (8a+\frac{1}{q_0}).
	\end{equation}
    Now, condition (\ref{enumitem:cond_f_4}) applied at $x_2$ gives us that
    \begin{align}
    	\frac{1}{q_0} &\leq f(x_2) - x_2 f'(x_2) \nonumber\\
    	&\leq f(x_2) - bx_2 - x_2\left(f'(x_2)-b\right) \nonumber\\
    	&\leq |f(x_2) - bx_2| + x_2 |f'(x_2)-b| \nonumber\\
    	\label{align:using_bound_g_prime}
    	&< a + (x_0-x_2) \left(8a+\frac{1}{q_0}\right) + x_2\left(8a + \frac{1}{q_0}\right)\\
    	&= a + x_0(8a+\frac{1}{q_0})\nonumber\\
    	\label{align:using_x_0_small}
    	&\leq 5a + \frac{1}{2q_0}
    \end{align}
    where (\ref{align:using_bound_g_prime}) follows from (\ref{eq:borne_f_x_2}) and (\ref{align:borne_g_prime}), and (\ref{align:using_x_0_small}) from the fact that $x_0$ is less than $\frac{1}{2}$. Hence, $\frac{1}{q_0} < 10a$. This contradicts the fact that $q_0$ satisfies $\frac{1}{q_0} \geq 10a$. As a conclusion, the bound $|f'(x) - b| \leq 28a$ must hold on $\left[0, \frac{1}{4}\right]$.
\end{proof}

We can now prove the main result of this section.

\begin{prop}
	\label{prop:potential_high_systolic_ratio_is_almost_linear}
     Let $J \in \mathcal{C}^\infty([-1; 1], \mathbb R)$ satisfying \ref{enumitem:cond_high_systolic_ratio_contact}, \ref{enumitem:cond_high_systolic_ratio_systole}, \ref{enumitem:cond_high_systolic_ratio_volume}. Then
     \[
     	|J(1)+J(-1)|< 224\int_{-1}^1 \left(J(t)-tJ'(t)\right) \mathrm{d}t.
     \]
\end{prop}

\begin{proof}
	First, integrating by part gives that
	\begin{align}
		\int_{-1}^1 (J(t)-tJ'(t)) \mathrm{d}t &= 2 \int_{-1}^1 J(t) \mathrm{d}t - \left(J(1)+J(-1)\right) \nonumber\\
		&=2\left(\int_{0}^{1}\left(J(t)-tJ(1)\right) \mathrm{d}t + \int_{-1}^{0}\left(J(t)+tJ(-1)\right) \mathrm{d}t\right).
		\label{eq:integrales_prop_high_systolic_linear}
	\end{align}
	Since $\left(\frac{J(x)}{x}\right)' = - \frac{J(x)-xJ'(x)}{x^2}$, condition \ref{enumitem:cond_high_systolic_ratio_contact} implies that the function $x \mapsto \frac{J(x)}{x}$ is decreasing. In particular, the functions $h^+ \colon x \mapsto J(x)-xJ(1)$ and $h^- \colon x \mapsto J(-x)-xJ(-1)$ are non-negative on $[0,1]$. Since both $h^+$ and $h^-$ satisfy $h^\pm(x) - x (h^\pm)'(x) > 0$, Lemma \ref{lemma:consequence_volume} gives the bound
	\[	
		\forall x \in (0,1), 0 < h^\pm(x) < 2\frac{\int_0^1h^\pm(t) \mathrm{d}t}{x}.
	\]
	It follows that
	\[
		\max \limits_{x \in \left[\frac{1}{4}, \frac{1}{2}\right]}h^\pm(x) < 8\int_{0}^1h^\pm(t)\mathrm{d}t.
	\]
	Letting $a = 8\left(\int_{0}^1(h^+(t)+h^-(t))\mathrm{d}t\right)$, and since both $h^+$ and $h^-$ are non-negative, we get
	\begin{equation}
		\max \limits_{x \in \left[\frac{1}{4}, \frac{1}{2}\right]} |J(x)-xJ(1)| < a
		\label{eq:borne_h+}
	\end{equation}
	and
	\begin{equation}
		\max \limits_{x \in \left[\frac{1}{4}, \frac{1}{2}\right]} |J(-x)-xJ(-1)| < a
		\label{eq:borne_h-}
	\end{equation}

	Moreover, we have that
	\begin{align}
		a &= 8\left(\int_{0}^1(h^+(t)+h^-(t))\mathrm{d}t\right) \nonumber \\
		\label{align:chgt_variable}
		&= 8\left(\int_{0}^{1}\left(J(t)-tJ(1)\right) \mathrm{d}t + \int_{-1}^{0}\left(J(t)+tJ(-1)\right) \mathrm{d}t\right) \\
		\label{align:egalite_integrales}
		&=4\left(\int_{-1}^1 (J(t)-tJ'(t)) \mathrm{d}t\right)
	\end{align}
	where (\ref{align:chgt_variable}) follows from the change of variable $t \in [0,1] \mapsto -t \in [-1,0]$ in the integral $\int_{0}^{1}h^-(t)\mathrm{d}t$, and (\ref{align:egalite_integrales}) from (\ref{eq:integrales_prop_high_systolic_linear}). Finally, \ref{enumitem:cond_high_systolic_ratio_volume} and (\ref{align:egalite_integrales}) give
	\begin{equation}
		\label{eq:utilisation_C3}
		a < \frac{1}{20}.
	\end{equation}

	The uniform bound (\ref{eq:borne_h+}) (resp. (\ref{eq:borne_h-})), together with (\ref{eq:utilisation_C3}), allows to apply Lemma \ref{lemma:consequence_systole} to $x \in [0,1] \mapsto J(x)$ (resp. $x \in [0,1] \mapsto J(-x)$) with $b = J(1)$ (resp. $b = J(-1)$). This gives the bounds
	\[
		|J'(x)-J(1)| < 28a
	\]
	and
	\[
		|J'(-x)+J(-1)| < 28a
	\]
	for every $x \in \left[0,\frac{1}{4}\right]$.
	In particular, we have
	\[
		|J'(0)-J(1)| < 28a
	\]
	and
	\[
		|J'(0)+J(-1)| < 28a.
	\]
	The triangular inequality and (\ref{align:egalite_integrales}) give the result.
\end{proof}

\section{Proof of Theorem \ref{thm:ineg_sys}}

\begin{proof}{(Theorem \ref{thm:ineg_sys})}
	The proof goes by contradiction: fix $\epsilon$ an arbitrarily small positive number, and assume there is a contact form $\alpha$ in $\Omega\left(g,(p_l,q_l)_{1 \leq l \leq L}\right)$ with volume less than $\epsilon$ and systole larger than one. In particular, every critical value $K_0$ of $K$ satisfies $|K_0| \geq 1$, since there is an integer $p$ such that $\frac{|K_0|}{p}$ is the minimal period of a closed Reeb orbit. This implies that for all $i \in I_0$, $[-1,1]$ is contained in $[K_i^-, K_i^+]$. Using Lemma \ref{lemma:volume_form_decomposition}, the following holds:
	\begin{align*}
		\vol(M \setminus K^{-1}([-1,1])) &\geq \int_{M \setminus (K^{-1}([-1,1]) \cup \mathcal A)} \sigma \wedge \pi^*\left(\tilde K \mathrm{d} \beta + \beta \wedge \mathrm{d} \tilde K\right) \\
		&=\int_{\tilde S \setminus \tilde K^{-1}([-1,1])}\left(\tilde K \mathrm{d} \beta + \beta \wedge \mathrm{d}\tilde K\right)\\
		&\geq \int_{\tilde S \setminus \tilde K^{-1}([-1,1])}\frac{\tilde K \mathrm{d} \beta + \beta \wedge \mathrm{d}\tilde K}{\tilde K^2}\\
		&=\int_{\tilde S \setminus \tilde K^{-1}([-1,1])}\mathrm{d}\left(\frac{\beta}{\tilde K}\right)\\
		&=\int_{\partial \left(\tilde S \setminus \tilde K^{-1}([-1,1])\right)} \frac{\beta}{\tilde K} \\
		&=\int_{\partial \tilde S} \frac{\beta}{\tilde K} + \sum_{i \in I_0} (J_i(1) + J_i(-1)).
	\end{align*}
	Moreover, $\int_{\tilde S \setminus \tilde K^{-1}([-1,1])}\frac{\tilde K \mathrm{d} \beta + \beta \wedge \mathrm{d}\tilde K}{\tilde K^2}$ is positive. Indeed, $\sigma \wedge \pi^* \left(\tilde K \mathrm{d} \beta + \beta \wedge \mathrm{d}\tilde K\right)$ coincides with the contact volume form on $\pi^{-1}(\tilde S)$. Hence, $\tilde K \mathrm{d} \beta + \beta \wedge \mathrm{d}\tilde K$ is a positive area form on $\tilde S$ given the quotient orientation, and so is $\frac{\tilde K \mathrm{d} \beta + \beta \wedge \mathrm{d}\tilde K}{\tilde K^2}$ on $\tilde S \setminus \tilde K^{-1}([-1,1])$. This gives
	\[
		\int_{\partial \tilde S} \frac{\beta}{\tilde K} + \sum_{i \in I_0} (J_i(1) + J_i(-1)) > 0.
	\]
	Lemma \ref{lemma:computing_euler}, together with the fact that $\alpha$ has volume less than $\epsilon$, gives
	\[
		0 < \sum_{i \in I_0} (J_i(1) + J_i(-1)) + e < \epsilon
	\]
	which becomes
	\begin{equation}
		\label{eq:volume_vs_euler}
		-e < \sum_{i \in I_0} (J_i(1) + J_i(-1))< -e+\epsilon.
	\end{equation}
	From now on, we assume that $\epsilon < \frac{1}{80}$. For each $i \in I_0$, we can apply Proposition \ref{prop:potential_high_systolic_ratio_is_almost_linear} to the restriction to $[-1,1]$ of $J_i$. Indeed,
	\begin{itemize}
		\item \ref{enumitem:cond_high_systolic_ratio_contact} follows from Corollary \ref{coro:identity_tau},
		\item \ref{enumitem:cond_high_systolic_ratio_systole} follows from Corollary \ref{coro:periodic_point} together with the condition on $\alpha$ to have systole larger than one,
		\item \ref{enumitem:cond_high_systolic_ratio_volume} follows from Lemma \ref{lemma:volume_J} and the assumption $\epsilon < \frac{1}{80}$.
	\end{itemize}

	This gives that for all $i \in I_0$, $|J_i(1) + J_i(-1)| < 224 \int_{-1}^1 J_i(x)-xJ_i'(x) \mathrm{d}x$. We get that
	\begin{align}
		|\sum_{i \in I_0}(J_i(1) + J_i(-1))| &< \sum_{i \in I_0}|J_i(1) + J_i(-1)| \nonumber\\
		&< 224\sum_{i \in I_0}\int_{-1}^1 \left(J_i(x)-xJ_i'(x)\right) \mathrm{d}x \nonumber\\
		\label{eq:sum_is_small}
		&< 224 \epsilon.
	\end{align}

	Hence, $\sum_{i \in I_0}(J_i(1) + J_i(-1))$ is in $]-e,-e+\epsilon[$ by (\ref{eq:volume_vs_euler}) and in $]-224 \epsilon, 224 \epsilon[$ by (\ref{eq:sum_is_small}).
	However those intervals do not intersect when $\epsilon$ is less than $\frac{|e|}{225}$. Since we assumed that $\epsilon$ was less than $\frac{1}{80}$, we get that there is no contact form $\alpha$ in $\Omega\left(g,(p_l,q_l)_{1 \leq l \leq L}\right)$ with volume less than $\min\left\{\frac{1}{80}, \frac{|e|}{225}\right\}$ and systole at least one. In other words, the systolic inequality
	\[
		\sys(\alpha)^2 \leq \max\left(80, \frac{225}{|e|}\right) \vol(\alpha)
	\]
	holds on $\Omega\left(g,(p_l,q_l)_{1 \leq l \leq L}\right)$. The result follows for $C = 225$.
\end{proof}

\bibliographystyle{plain}
\bibliography{biblio}
\end{document}